\documentclass{amsart}
\usepackage{amsfonts,amscd}
\newtheorem{claim}{}[section]
\newtheorem{theorem}[claim]{Theorem}
\newtheorem{lemma}[claim]{Lemma}
\newtheorem{proposition}[claim]{Proposition}
\newtheorem{corollary}[claim]{Corollary}

\newtheorem{definition}[claim]{Definition}

\def\proclaim #1. #2\par{\medbreak
\noindent{\bf#1.\enspace}{\sl#2}\par\medbreak}
\makeatother

\DeclareMathOperator{\ca}{C^*-\text{algebra}}

\DeclareMathOperator{\oas}{\text{operator algebras}}

\DeclareMathOperator{\auoa}{\text{approximately unital operator algebra}}
\DeclareMathOperator{\auoas}{\text{approximately unital operator algebras}}
\DeclareMathOperator{\cas}{C^*-\text{algebras}}

\begin{document}
\title[Logmodularity and isometries of
operator algebras] {Logmodularity and isometries of operator
algebras}

\author{David P. Blecher}
\address{Department of Math, University of Houston, Houston, TX77204-3476}
\email{dblecher@math.uh.edu}

\author{Louis E. Labuschagne}
\address{Department of Math, Applied Math, and Astronomy, P.O. Box
392, 0003 UNISA, South Africa}

\email{labusle@unisa.ac.za}

\thanks{April 30, 2002.}

\thanks{This research was supported in part by grants from
the National Science Foundation and the University of South
Africa.}

\vspace{30 mm}

\begin{abstract}
We generalize some facts about function algebras
to operator algebras, using the `noncommutative 
Shilov boundary' or $C^*$-envelope first considered 
by Arveson.  In the first part we 
study and characterize complete isometries
between operator algebras.  In the second part
we introduce and 
study a notion of logmodularity for operator algebras.
We also give a result on conditional expectations.
Many miscellaneous applications are given.   \end{abstract}

\maketitle

\section{Introduction.}

The main topic of our paper is the study of linear maps $T : A
\rightarrow B$ between operator algebras. By an `operator algebra'
we mean a uniformly closed algebra of operators on a Hilbert space
$H$, and we usually assume that such algebras are unital (that is,
contain an identity of norm 1), or are  {\em  approximately
unital} (that is, contain a contractive two-sided approximate
identity (c.a.i.)).     Operator algebras may be viewed as
`noncommutative function algebras'. Indeed every function algebra
or uniform algebra  on a compact space $K$ containing constant
functions, is a closed unital subalgebra of the commutative $\ca$
$C(K)$, and hence is a unital operator algebra.

It is nowadays commonly recognized that to study a general
subalgebra $A \subset B(H)$ it is necessary not only to take into
account the norm, but also the natural norm on the matrix spaces
$M_n(A) \subset M_n(B(H)) \cong B(H^{(n)})$.  This has been one of
the key perspectives of operator space theory, since Arveson's
pioneering work; and the present paper may in some sense
be regarded as an extended series of observations proceeding from
Arveson's papers \cite{Arv1,Arv2,Arvan}.

Hence we are not interested
here in bounded linear transformations, rather we look
for the completely bounded maps - where the adjective
`completely' means that we are applying our maps to
matrices too.  Thus if $T : X \rightarrow Y$, then $T$
 is completely contractive if and only if
 $$\Vert [T(x_{ij})] \Vert
\leq \Vert  [x_{ij}] \Vert $$
for all $n \in \mathbb{N}$ and
$[x_{ij}] \in M_n(X)$.
 We say that $T$ is
{\em completely isometric}
if $\Vert [T(x_{ij})] \Vert = \Vert [x_{ij}] \Vert$
for all $n \in \mathbb{N}$ and
$[x_{ij}] \in M_n(X)$.

In much of this paper we are interested in
complete isometries between operator algebras, and in 
generalizing facts for isometries between function algebras 
(note that isometries between function algebras 
are completely isometric (and vice versa),
as may be seen from \cite{P} Theorem 3.8). 
We now briefly review some of these facts.
The classical Banach-Stone theorem characterizes onto linear
isometries between $C(K)$ spaces, that is, between commutative
 unital $\cas$.   There is an analogous result for function
algebras (see e.g. p. 144 in \cite{Ho}; \cite{Nag}).   These
results have been extended in many directions (see \cite{JP} for a
 survey).
The not necessarily surjective isometries
between $C(K)$ spaces were characterized by
 Holsztynski, whereas in the case of function algebras there is the
 following more general result from \cite{Math} (for closely related
 results see \cite{Nov} and particularly \cite{AF} Section 6):

\begin{theorem} \label{mat} \cite{Math}
A linear map $T :  A \rightarrow B$ between
uniform algebras, is an isometry if and only if $T$ is
contractive, and there exists a closed subset $E$ of $\partial B$,
and two continuous functions $\gamma : E \rightarrow \mathbb{T}$
and $\varphi : E  \rightarrow \partial A$, with $\varphi$
surjective, such that
$$T(f)(y) \; = \; \gamma(y) f(\varphi(y))$$ for all $y \in E$.
\end{theorem}

Here $\mathbb{T}$ is the unit circle, and
$\partial B$ is the Shilov boundary of $B$ (which equals
$K$ if $B = C(K)$).
Informally, the result is saying
that for any isometry $T$
there is a certain `part $E$ of the space $B$ acts on',
such that $T$ restricted to this part has a particularly nice form,
a form which amongst other things ensures that $T$ is an isometry.
The action of $T$ on the `complementary part' plays little
role in the fact that $T$ is an isometry.
If $T$ is {\em unital} (that is, $T(1) = 1$), 
then the $\gamma$ may be omitted in the
theorem, and then $T$ compressed to $E$ is the composition
operator $f \mapsto f \circ \varphi$, which is
an isometric homomorphism.

For surjective isometries between
$\cas$ the first `noncommutative Banach-Stone theorem'
is due to R. V. Kadison \cite{Ka}.    Related isometric results 
were obtained for example by Harris and Arazy and Solel (see 
\cite{AS}). In  Part I of our paper we consider 
unital or approximately unital  operator algebras $A$ and $B$,
and (not necessarily surjective) linear
maps $T : A \rightarrow B$, and we will establish several
criteria which are each equivalent to $T$ being a
linear complete isometry.   These are generalizations
of \ref{mat} because as we said earlier,
 any function algebra
is a unital operator algebra; and moreover
isometries between function algebras are completely isometric.
In the result above we saw the importance of the
Shilov boundary; in the noncommutative case we will need to use
the {\em noncommutative Shilov boundary} or
{\em $C^*$-envelope} $C^*_e(A)$ of an operator algebra $A$.   This
will be described in more detail later towards the end of Section 1.

One such characterization proceeds as follows. Suppose that $A, B$
are operator algebras, and suppose that $B$
is a subalgebra of a $\ca$ $C$.
Suppose that  $T : A \rightarrow B$ is a linear map.
For the purposes of additional clarity here
we assume that $A, B$ and $T$ are unital
(in the
general case one needs to add a partial isometry, and a
second projection $q$, to the
statement that follows).
Then $T$ is a complete isometry if and only if there is a
 projection $p$
in $B(H_u)$ where $H_u$ is the Hilbert space of the universal
representation of $C$, such that firstly,
$$p T(\cdot) = T(\cdot) p = p T(\cdot) p = \pi(\cdot)$$ for
a completely isometric
homomorphism $\pi$ which is the restriction to $A$ of a
1-1 *-representation $C^*_e(A) \rightarrow B(p H_u)$;
and secondly,
the map $S = (1-p) T(\cdot) (1-p)$ is a complete contraction.
This is saying that up to a choice of orthonormal bases,
$T$ is of `block diagonal' form
$$\left[ \begin{array}{ccl}
\pi(\cdot) & 0 \\ 0 & S(\cdot) \end{array} \right].$$

In Part I we also consider particularly interesting
subclasses of complete isometries between  $\oas$.
In particular, we investigate a class which was introduced in the
commutative case by Matheson under the name {\em type 1 isometry}.

As we have said, an important tool
in the noncommutative case is the
$C^*$-envelope $C^*_e(A)$ of an operator algebra $A$.   However there are two
basic obstacles associated with this tool, which we address (with very partial
success) in the second half of our paper.   The first obstacle is that of
{\em identifying} $C^*_e(A)$ for a given   operator algebra
$A$.   Following the classical case, we introduce a notion
of {\em logmodularity} for operator algebras, and
show that for a logmodular subalgebra $A$ of a $\ca$ $B$,
the $C^*$-envelope  of $A$ is $B$.
We develop some basic theory of logmodular algebras, giving for
example a result
which is a noncommutative version of the uniqueness of
representing measures.  Indeed this result ensures that
a completely contractive homomorphism $\theta$ defined on a
logmodular operator algebra has at most one
completely contractive extension to $C^*_e(A)$.    This
extension will be completely isometric if $\theta$ is.

A good example of a logmodular algebra in the commutative case is
the Hardy space $H^\infty(\mathbb{D})$; and in the noncommutative
case the {\em noncommutative $H^\infty$ algebras} introduced by
Arveson (also known as finite maximal subdiagonal algebras).  In
section 5 we
make a few observations about these algebras using the facts
in the previous paragraph.  In this section we also
prove some general results on conditional expectations of
non-selfadjoint operator algebras.

The second obstacle associated with the $C^*$-envelope arises when
 $C^*_e(A)$ is rather complicated.   Again a good example of this
 is $H^\infty(\mathbb{D})$.   In the light of
the known characterizations of isometries of other common
 function algebras in terms of {\em composition operators}
 (see \cite{SM} for example), one would hope for a
 characterization in which analytic self-maps of the
 open disk play a major role (as opposed to
 self-maps of largely unhelpful Shilov boundary).
 However for this purpose our results on
 complete isometries from Part I are often not helpful.  To
 illustrate this point we recall that there is no such known
characterization of isometries $H^\infty(\mathbb{D})
\rightarrow H^\infty(\mathbb{D})$.   
In section 3 we make a few observations concerning
isometries of $H^\infty(\mathbb{D})$.
Our main interest however is in the noncommutative
versions, and we are currently working with A. Matheson
in this direction.

\vspace{4 mm}

We end this introduction with some notation and basic facts. For
more details on operator spaces we refer the reader to
\cite{ERbook,P,Pisbook}. We write $H, K, L$ for Hilbert spaces.
Perhaps confusingly, we will use the symbols $I, J, K$ for ideals
in a $\ca$.   All ideals are taken to be uniformly closed. A
projection on a Hilbert space, or in a $C^*$-algebra, will mean an
orthogonal projection.  Otherwise we use the word projection for a
linear self-map $P$  with $P \circ P = P$. If $Y$ is a subspace of
$X$ we write $q_Y$ or $q$ for the natural quotient map $X
\rightarrow X/Y$.

For a set of operators on a Hilbert space, or for a subset of an
operator algebra we write $S^\ast$ for the subset corresponding to
the set of adjoints of elements of $S$; $S_+$ to be the set of
$x \in S$ with $x \geq 0$; and $S^{-1}$ for the invertible
elements of $S$ whose inverse is also in $S$.     We write $|b|$
for $(b^* b)^{\frac{1}{2}}$.  To avoid possible confusion the
Banach dual of a Banach space $X$ will be denoted by $X^\star$.
The canonical map $X \rightarrow X^{\star \star}$ is written as
$\; \hat{} \;$.

If $X, Y$ are subspaces of a Banach algebra, we write $XY$ for the
{\em uniform closure } of the set of finite sums of products of
the form $x y$ for $x \in X, y \in Y$.

An  {\em operator system} is a unital selfadjoint subspace $X$ of
a unital $\ca$ or of $B(H)$ (i.e. $X$ contains the identity
element, and $x \in X$ if and only if $x^* \in X$). A {\em unital
operator space} has the same definition except we do not require
that $X^* = X$. The appropriate morphisms between operator systems
are unital completely positive maps.
Rather than state the (obvious) definition
of a completely positive map, we will simply use the fact
that a unital linear map $S$ between operator systems is
completely contractive if and only if it is completely positive,
in which case it is *-linear: i.e.  $S(v^*) = S(v)^*$. See
\cite{Arv1,P} for proofs.

\begin{lemma}  \label{Arv} (Arveson, see e.g. \cite{Con}) $\;$
Suppose that $V$ is a subspace of an operator system which
contains the identity element, and suppose that $\psi : V
\rightarrow B(K)$ is a unital contraction with range $W$. Then
there exists a unique positive map between the  operator systems
$V + V^*$ and $W + W^*$, which extends $\psi$. If $\psi$ is a
unital complete contraction (resp. unital complete isometry), 
then the
unique positive extension to $V + V^*$ is a completely positive
map (resp. a unique complete order isomorphism) between the
operator systems $V + V^*$ and $W + W^*$.
\end{lemma}

As an aside we note that in the context of $H^p$-spaces the above
result of Arveson seems to translate to a principle whereby
composition operators on $H^p(\mathbb{D})$ spaces may be extended
to positive integral operators on the enveloping
$L^p(\mathbb{T})$-spaces. (See \cite{Sar}.)

Now if $T  : X \rightarrow A$, and $S : X \rightarrow B$, are
unital complete isometries of an abstract unital operator space
$X$ into $C^*$-algebras $A$ and $B$, then by the above result
there exists a unique complete order isomorphism between the
operator systems $T(X) + T(X)^*$ and $S(X) + S(X)^*$ which extends
the map $S \circ T^{-1}$ from $T(X) \rightarrow S(X)$.
As a corollary of this, one may easily check that if $A$ is
an abstract unital operator algebra, then the two spaces
$\Delta(A) = A \cap A^*$ (the diagonal of $A$), and $A + A^*$, are
`well defined' independently of the particular $C^*$-algebra
containing $A$ as a unital subalgebra.   Note that $\Delta(A)$ is
a $C^*$-algebra (indeed is a $W^*$-algebra if $A$ is a `weak*
closed operator algebra'), and $A + A^*$ is an operator system.
See also \cite{AS} for some interesting related facts.

If $X$ is a unital operator space then there exists a {\em
$C^*$-envelope} of $X$, namely a pair $(B,j)$ consisting of a
unital $\ca$ $B$ and a unital complete isometry $j : X \rightarrow
B$ whose range generates $B$ as a $\ca$, with the following
universal property:  For any other pair $(A,i)$ consisting of a
unital $\ca$ and a unital complete isometry $i : X \rightarrow A$
whose range generates $A$ as a $\ca$, there exists a (necessarily
unique, unital, and surjective) *-homomorphism $\pi : A
\rightarrow B$ such that $\pi \circ i = j$. This we call the
Arveson-Hamana theorem \cite{Arv1,Ham}, and as customary we write
$C^*_e(X)$ for $B$ or $(B,j)$ (it is essentially unique, by the
universal property).

If $A$ is a unital operator algebra, then it is easy to see that
$j$ in the last paragraph is forced to be a completely isometric
homomorphism. Thus we consider $A$ to be  a unital subalgebra of
$C^*_e(A)$.  If $A$ is an $\auoa$ then one may define $C^*_e(A)$
to be the $\ca$ generated by $A$ inside the $C^*$-envelope of the
unitization of $A$.  Then one may easily see that $C^*_e(A)$ has
the desired universal property.

For the purposes of this paper we define a {\em triple system} to
be a (uniformly closed) subspace $X$ of a $\ca$ such that $X X^* X
\subset X$.
The important structure on a triple system is the `triple product'
$x y^* z$.  A `triple subsystem' is a uniformly closed vector
subspace of a triple system which is closed under this triple
product.

It is well known  that $X X^* X = X$ for a triple system $X$.
Also, it is clear that $X X^*$ and $X^* X$ are $C^*$-algebras,
which we will call the left and right $C^*$-algebras of $X$
respectively, and $X$ is a $(X X^*) - (X^* X)$-bimodule. A linear
map $T : X \rightarrow Y$ between triple systems is a {\em triple
morphism} if $T(x y^* z) = T(x) T(y)^* T(z)$ for all $x,y,z \in
X$.   Triple systems are operator spaces, and
 triple morphisms behave very similarly to *-homomorphisms
between $\cas$: triple morphisms are automatically completely
contractive and have closed range. A triple morphism is completely
isometric if it is 1-1.  The kernel of a triple morphism on $X$ is
a `triple ideal' (that is, a uniformly closed
 $(X X^*)-(X^* X)$-subbimodule).  The quotient of a
triple system by a triple ideal is a triple system in an obvious
way.  If one factors a triple morphism by its kernel one obtains a
1-1 triple morphism on the quotient triple system.

If $X$ is any operator space, then there exists a {\em triple
envelope} of $X$, namely a pair $(Z,j)$ consisting of a triple
system $Z$ and a linear complete isometry $j : X \rightarrow Z$
whose range generates $Z$ as a triple system (that is there exists
no nontrivial closed triple subsystem containing $j(X)$), with the
following universal property:  For any other pair $(W,i)$
consisting of a triple system and a  complete isometry $i : X
\rightarrow W$ whose range generates $W$ as a triple system, there
exists a (necessarily unique and necessarily surjective) triple
morphism $\pi : W \rightarrow Z$ such that $\pi \circ i = j$.

This theorem dating to the `80's is due to Hamana \cite{Ham3}, and
we write ${\mathcal T}(X)$ for $Z$ or $(Z,j)$ (as before it is
essentially unique, by the universal property).   The 
triple envelope or $C^*$-envelope is usually defined as a subspace
of the {\em injective envelope} $(I(X),j)$.  We will also need 
this latter  envelope, but since our 
use of it is quite limited we will not take the time to 
define it here, and instead refer the reader to \cite{Ham3} 
or \cite{ERbook} for more details 
(the forthcoming
revision of \cite{P} focuses extensively on this
topic too).   See also
\cite{BShi} for a more thorough discussion of the  
concepts in the last several paragraphs.

\vspace{5 mm}

\begin{center}
{\Large Part I.   Characterizations of complete isometries.}
\end{center}

\vspace{1 mm}

\section{General Banach-Stone theorems between nonselfadjoint
operator algebras}

Suppose that $T : A \rightarrow B$ is a surjective unital complete
isometry between unital operator algebras $A$ and $B$. Then by the
universal property of the $C^*$-envelope $T$ extends to a unital
complete isometry $\tilde{T}$ between the $C^*$-envelopes
$C^*_e(A)$ and $C^*_e(B)$. 
 By the universal 
property $\tilde{T}$ is a *-homomorphism
(or this may be seen by a more simple Banach-Stone theorem
for $C^*$-algebras such as 5.2.3 in \cite{ERbook}).
Consequently $T$ is also a homomorphism. In
fact even if $T$ is not unital the same argument works (modulo a
few technical details) to show that:

\begin{theorem} \label{bstoa}
(\cite{Arv1,Arv2,Bcomm,ERnsa} and B.1 in \cite{BShi}) Suppose that
$T : A \rightarrow B$ is a linear complete isometry between
approximately unital operator algebras.    Then $T$ is a
surjective complete isometry if and only if $T = u \theta(\cdot)$,
where $u$ is a unitary in the diagonal $\Delta(M(B))$ of the
multiplier algebra of $B$, and where  $\theta : A \rightarrow B$
is a surjective completely isometric homomorphism which is a
restriction of a *-isomorphism between the $C^*$-envelopes of $A$
and $B$.  Furthermore, the surjective complete isometry $T$ is
unital if and only if the $u$ above equals $1_B$, and if and only
if $T$ is a homomorphism.
\end{theorem}

In the remainder of this section, and in the following section, we
attempt to find analogous results for nonsurjective complete
isometries between nonselfadjoint operator algebras; and give some
applications.  These will generalize Theorem \ref{mat} stated in
the introduction.

The first set of characterizations we will collect in the
following theorem, which is a generalization of Theorem 3.1 in
\cite{BH}.  We apologize in advance to the reader - since this is
quite similar to the analogous result for selfadjoint operator
algebras in \cite{BH}, we cannot justify repeating all the ideas
and notation established there.    Thus we must  refer the reader
to this paper for further details regarding our notation and
proof.

Suppose that $B$ is a subalgebra of a $\ca$ $C$. A reducing
$\wedge$-compression of a map $T : X \rightarrow B$ is a map $R :
X \rightarrow C^{\star \star}$ such that there exist projections
$e, f \in C^{\star \star}$ such that $$R = \widehat{T(\cdot)} e =
f \widehat{T(\cdot)} = f \widehat{T(\cdot)} e.$$ Thus if $C$ is
represented on its universal Hilbert space $H_u$, then there are
subspaces $L, M$ of $H_u$ (corresponding to $p, q$) such that $R =
P_M T(\cdot)_L$, and $T = R \oplus S$ for a completely contractive
$S : X \rightarrow B(L^\perp,M^\perp)$.

\begin{theorem}  \label{ncmat}
Let $T :  A \rightarrow B$ be a linear completely contractive map
between approximately unital operator algebras. Suppose that $B$
is a subalgebra of a $\ca$ $C$ (for example $C = C^*_e(B)$).  Then
the following are equivalent:
\begin{itemize}
\item [(i)]  $T$ is a complete isometry,
\item [(ii)]  there is a left ideal $J$, and a right ideal $K$, of
$C$, such that  $q_{J + K} \circ T$ is the restriction to $A$ of a
1-1 partial triple morphism; that is the restriction of a linear
map $S: C^*_e(A) \rightarrow C/(J+K)$ which on composition with
the canonical injection $C/(J+K) \rightarrow (1 - q)C^{\star
\star}(1-p)$ is a triple morphism. (Here $q$ and $p$ are the
so-called support projections of $J$ and $K$ respectively. See
\cite{BH} for the details.)
\item [(iii)]
$T$ has a reducing $\wedge$-compression $R : A \rightarrow
C^{\star \star}$ which is the restriction of a 1-1 triple morphism
from $C^*_e(A)$ into $C^{\star \star}$.
\item [(iv)] there is a projection
$p$ in $C^{\star \star}$, a  partial isometry $u$ in $C^{\star
\star} (1-p)$, and a 1-1 *-homomorphism $\theta : C^*_e(A)
\rightarrow (1-p) C^{\star \star} (1-p)$, such that
$$\widehat{T(a)} (1-p) = u \theta(a),$$ and $u^* u \theta(a) =
\theta(a)$, for all $a \in A$.
\end{itemize}

If $A$ is unital these are equivalent to
\begin{itemize}
\item [(v)] there exists a closed left
ideal $J$ of $C$, such that for all $a \in A$ we have $$q_J(T(a))
\; = \; u  \pi(a) \quad \mbox{and} \quad \pi(a) \; = \; u^*
  q_J(T(a)) $$
where $q_J : C \rightarrow C/J$ is the canonical quotient map,
$\pi : A \rightarrow C/(J + J^*)$ is a completely isometric
homomorphism which is the restriction of a 1-1 *-homomorphism
$C^*_e(A) \rightarrow M_{00}(C/(J + J^*))$ (see \cite{BH} for
definitions), and $u$ is a partial isometry in $B/J \subset C/J$
which may be taken to be $u = q_J(T(1))$,
\item [(vi)]  there exists a $C^*$-subalgebra $D$ of $C$
containing $T(A)^* T(A)$, a closed two-sided ideal $I$ in $D$, and
a surjective *-isomorphism $\pi : C^*_e(A) \rightarrow D/I$, such
that $$q_N(T(a)) \; = \; U \; \; \pi(a)$$ for all $a \in A$, where
$q_N$ is the canonical quotient triple morphism from the triple
subsystem $Z = T(A) D$ of $C$, onto the quotient of $Z$ by the
triple ideal $N = Z I$, and $U$ is a unitary in the triple system
$Z/N$.  One may take $U = q_N(T(1))$.
\end{itemize}
 \end{theorem}

\begin{proof}   This is very similar to the proof of
3.1 in \cite{BH}, so we will be brief and omit 
some of the proof details.  In particular, it is a
simple exercise that any of (ii)-(vi) implies (i). Let $Z$ be the
triple system in $C$ generated by $T(A)$.   In \cite{BShi} section
4 it is shown that $C^*_e(A)$ is a `triple envelope' of $A$.  Thus
by the universal property of this envelope mentioned in the
introduction, there exists a surjective triple morphism $\rho : Z
\rightarrow C^*_e(A)$ such that $\rho(T(a)) = a$ for all $a \in
A$.   Set $N = \mathrm{Ker}\, \rho$, and obtain a surjective
triple isomorphism $\xi : C^*_e(A) \rightarrow Z/N$.   Appealing
to Lemma 2.10 in \cite{BH} almost immediately gives (vi). We may
then follow the proof of (i) $\Rightarrow$ (iii) of 3.1 in
\cite{BH}, replacing $B$ there by $C$, and so on, to obtain a
completely isometric partial triple morphism $\theta$ of
$C^*_e(A)$ into $C/(J + K)$. The restriction of $\theta$ to $A$ is
$q_{J + K} \circ T$.   This gives (ii). We also obtain as in 3.1
of \cite{BH}, a 1-1 triple morphism $\tilde{T} : C^*_e(A)
\rightarrow C^{\star \star} (1-p)$ taking $a \in A$ to
$\widehat{T(a)}(1-p) = (1-q) \widehat{T(a)}(1-p) = (1-q)
\widehat{T(a)}$.   This gives (iii). Applying 2.10 of \cite{BH}
gives $\widehat{T(a)}(1-p) = u \theta(a)$ for a *-homomorphism
$\theta : C^*_e(A) \rightarrow (1-p) C^{\star \star} (1-p)$
thereby yielding (iv).
The proof of (v) is then similar to the proof that (i) implies
(vi) in 3.1 of \cite{BH}.
 \end{proof}

It will be readily seen how \ref{mat} follows immediately from for
example (v) of our theorem.   In particular as stated in the
introduction, under the hypotheses of \ref{mat}, $T$ is completely
isometric and we may apply \ref{ncmat}.   Since we find ourselves
in a commutative context in 1.1, the (now two-sided) ideal $J = J
+ J^*$ corresponds to a set of continuous functions vanishing on
some closed subset $E$ of $\partial B$.  Thus $C/(J + J^*) =
C(\partial B)/J$ is then a copy of $C(E)$. The 1-1
$\ast$-homomorphism from $C(\partial A)$ into $M_{00}(C/(J + J^*))
= C(E)$ is then induced by a continuous surjection $\varphi : E
\rightarrow \partial A$, with the partial isometry corresponding
to a continuous $\gamma : E \rightarrow \{0 \} \cup \mathbb{T}$.
Therefore (v) of the above theorem translates to the claim that
$$T(f)(y) = \gamma(y) f(\varphi(y)) \quad \mbox{and} \quad
f(\varphi(y)) \; = \; \overline{\gamma(y)} T(f)(y)$$ for all $y
\in E$ and $f \in A$. It follows from the second equation that
$\gamma : E \rightarrow \mathbb{T}$.

In the proof of (iv)
in \ref{ncmat}, we obtained projections $p, q$ in $C^{\star
\star}$ such that $$\widehat{T(\cdot)}(1-p) = (1-q)
\widehat{T(\cdot)} = (1-q) \widehat{T(\cdot)} (1-p)$$ and saw that
this expression is the restriction to $A$ of a 1-1 triple morphism
$\tau : C^*_e(A) \rightarrow (1-q) C^{\star \star}(1-p)$. From
this it is clear that $$T(ab)(1-p) = T(a) T(1)^* T(b)(1-p)\; \; \;
\; \; \; \; \; (*)$$ for all $a,b \in A$.   We will use this fact
later on.

\begin{corollary} \label{ncmatco}
 Let $T : A \rightarrow B$ be a unital completely contractive
linear map between unital operator algebras.  Suppose that $B$
 is a unital subalgebra of a unital $\ca$ $C$. Then saying that $T$ is
 a complete isometry is equivalent to any one of (ii)-(vi) in the
 previous theorem, but with the following changes: One may omit all
 mention of $u$, change `triple morphism' to `*-homomorphism', and
 $K$ to $J^*$, in (ii), and note that the expression
 $\widehat{T(a)} (1-p)$ in (ii) and (iv) now coincides with $(1-p)
 \widehat{T(a)} (1-p)$. One may change
(vi) to read:
 there exists a $C^*$-subalgebra $D$ of $C$
containing $T(A)$,
a closed two-sided ideal $I$ in $D$, and a
 surjective *-isomorphism $\pi : C^*_e(A) \rightarrow D/I$, such that
$$q_I(T(a)) \; = \; \pi(a)$$
for all $a \in A$.
\end{corollary}

The appropriate parts of the last theorem 
and corollary may also be stated in terms of block diagonal
matrices as in the introduction:   Suppose for example that $T : A
\rightarrow B$ is a unital complete isometry between unital
operator algebras, and suppose that $B$ is a unital subalgebra of
a unital $\ca$ $C$.  Let $H_u$ be the Hilbert space of the
universal representation of $C$, thus we have *-subalgebras $C
\subset C^{\star \star} \subset B(H_u)$.   Then corresponding to
the projection $p$ there is a subspace $K$ of $H_u$ such that $P_K
T(\cdot) = T(\cdot) P_K = P_K T(\cdot) P_K$, and if we define $R :
A \rightarrow B(K)$ by these equal expressions, then $R$ is a
completely isometric homomorphism (and $R$ is also  the
restriction of a 1-1 *-homomorphism defined on $C^*_e(A)$).   Then
$T = R \oplus S$ for a completely contractive unital map $S$ as in
the introduction.

A special case of the above worthy of consideration follows (For
the proof one need only note that the simplicity requirement
ensures that $I = \{0\}$ in \ref{ncmatco}.):

\begin{corollary}  \label{eaco}
Let $T : A \rightarrow B$ be a unital complete isometry between
unital operator algebras $A$ and $B$.  Suppose that $B$ is a
unital subalgebra of a $C^*$-algebra $C$, and suppose further that
the $C^*$-subalgebra of $C$ generated by $T(A)$ is simple (for
example if $T$ is a map into $M_n$ whose range generates $M_n$).
Then $T$ is automatically a homomorphism, indeed $T$ is the
restriction of a 1-1 *-homomorphism $C^*_e(A) \rightarrow C$.
\end{corollary}

We proceed with some typical applications of \ref{ncmat}:

\begin{corollary}  \label{uptr}  Let $T_n$ be the upper triangular
$n \times n$ matrices.  Then a linear map $\varphi : T_n
\rightarrow M_m$ is a complete isometry  if and only if there
exist unitary matrices $U$ and $V$ of appropriate sizes, and a
completely contractive map $S : T_n \rightarrow M_{m-n}$, such
that $\varphi(A) = U \mathrm{diag} \{ A , S(A) \} V$ for all $A
\in T_n$.
\end{corollary}

\begin{proof}  By the proof of
Theorem \ref{ncmat} there exist projections $p, q$ in $M_m$, and a
1-1 triple morphism $\theta : M_n  = C^*_e(T_n) \rightarrow (1-q)
M_m (1-p)$  such that $$(1-q) \varphi(A) = \varphi(A) (1-p) =
(1-q) \varphi(A) (1-p) = \theta(A)$$ for all $A \in T_n$.   The
result then follows from Lemma 3.5 of \cite{BH} (as in 3.6 of
\cite{BH}).
\end{proof}

\begin{corollary}  Let $A$ be a unital operator algebra
which is linearly completely isometric to a subspace of $M_n$, for
a finite integer $n$.   Then $A$ is completely isometrically
homomorphic to a subalgebra of $M_n$, and hence $A$ is completely
isometrically homomorphic to a unital subalgebra of $M_m$, for
some $m \leq n$.
\end{corollary}

\begin{proof}   This follows immediately from \ref{ncmat}
(iv), since (iv) gives a completely isometric homomorphism $A
\rightarrow M_n$, from which the result is clear. A more direct
proof may also be given.
\end{proof}

The ideas in the last few results give a practical way of
computing the $C^*$-envelope (which in this case equals the
injective envelope) of a
unital subalgebra $A$ of $M_n$.  Namely,
first one replaces $M_n$ by the unital $*$-subalgebra $B$ generated by
$A$.  By the abstract principles used above, there is
a central projection $p$ in $B$ with $B (1-p) = C^*_e(A)$.
Since $B$ is unitarily equivalent to (up to multiplicity) a direct sum 
$\oplus_k \; M_{n_k}$ of full matrix 
algebras $M_{n_k}$, the central projection $p$ is 
unitarily equivalent to a direct sum of identity matrices $I_{n_i}$ or
$0_{n_i}$ matrices. 
A similar idea works to compute the triple or injective envelope of
a subspace of $M_n$, except that
now we have two projections $p$ and $q$.
It is instructive, and is also a good source of counterexamples,
to select a set of two or three matrices in 
$M_n$, and then to compute the $C^*$-envelope or triple 
envelope of their span.
  
We close this section with a
related fact clarifying the relationship between
the categories of unital operator spaces and
$\auoa$:

\begin{corollary} \label{auint}  If $A$ is an $\auoa$,
and if $A$ is linearly completely isometric to a unital operator
space, then $A$ is a unital algebra.
\end{corollary}

\begin{proof}
If $V$ is the unital operator space, then it follows that the
triple envelope of $A$ is triple isomorphic to the triple envelope
of $V$.  The latter may be taken to be $C^*_e(V)$, which  is
unital \cite{Ham}. Also in \cite{BShi} we showed that the triple
envelope of an approximately unital operator algebra  $A$ may be
taken to be $C^*_e(A)$ 
(part of this was observed by C. Zhang too).  
It is now easy to see that $C^*_e(A)$ is
unital, whence $A$ is also unital. (The triple morphism from
$C^*_e(V)$ onto $C^*_e(A)$ takes $V$ onto $A$ and $1 \in V$ onto
$1 \in A$.)
\end{proof}

\section{Characterizations of particular classes of isometries}

If one considers Theorem \ref{mat}, a `nicest class' of isometries
suggests itself: namely those for which $E = \partial A$.   A
moment's thought (using Lemma 2.10 of \cite{BH}) shows that this
class coincides with the class of maps $T : A \rightarrow B$ which
are the restrictions of a 1-1 triple morphism $C(\partial A)
\rightarrow C(\partial B)$. We shall call these the {\em Shilov
isometries}.  More generally, if $A$ and $B$ are $\auoas$ we say
that a map $T : A \rightarrow B$ is a {\em Shilov isometry},
 if $T$ is the restriction to $A$
of a 1-1 triple morphism $C^*_e(A) \rightarrow C^*_e(B)$.

The following two results, which are fairly superficial,
we state to afford a comparison with subsequent results.

\begin{lemma} \label{csca}  Let $T :  A \rightarrow B$
be a linear map between unital $\cas$.   The following are
equivalent:
\begin{itemize}
\item [(i)]  $T$ is a Shilov isometry (that is, a
1-1 triple morphism),
\item [(ii)]  $T = u \pi(\cdot)$ for a 1-1 *-homomorphism
$\pi : A \rightarrow B$, and a partial isometry $u \in B$ such
that $u^* u = \pi(1_A)$,
\item [(iii)]  $T$ is a complete isometry
with $\{ 0 \} = J \subset C = B$ in Theorem \ref{ncmat} (v).
\end{itemize}
We may take $u = T(1_A)$. If these equivalent conditions hold,
then $T(A)$ contains a unitary (resp. isometry, coisometry) of $B$
if and only if $u = T(1_A)$ is a unitary (resp. isometry,
coisometry) of $B$. Indeed if $1_B \in T(A)$ then $T(1_A)$ is a
unitary of $B$ and $T(A) = \pi(A)$, so that Ran $T$ is a
$C^*$-subalgebra of $B$, and also $\pi$ is unital. If in addition
$T(1_A) = 1_B$, then $T$ is a *-homomorphism.
\end{lemma}

\begin{proof}   The equivalence of (i)-(iii) we leave as
an exercise (using Lemma 2.10 in \cite{BH}). Set $u = T(1)$. We
will use the fact that if $R, S$ are contractions between Hilbert
spaces with $R S = I$, then $S$ is an isometry and $R = S^*$. Thus
if $T(A)$ contains a coisometry $v = T(a_0)$, then $v = u
\pi(a_0)$, giving $1_B = u (\pi(a_0) v^*)$. By the fact  above
about Hilbert space operators, $u$ is a coisometry.  By symmetry
we obtain the assertions for isometries and unitaries.

If $1_B \in T(A)$, then $u$ is a unitary of $B$ by the above.
Since $T$ is  a triple morphism, it is evident that Ran $T$ is a
$C^*$-subalgebra of $B$.  Write $1_B = T(a_0)$ for an $a_0 \in A$.
Then $\pi(a_0) = u^* u \pi(a_0) = u^* T(a_0) = u^*$, so that $u^*$
and $u$ are in the $\ca$ $\pi(A)$. Hence $T(A) = u \pi(A) =
\pi(A)$. Also $\pi(1_A) = u^* T(1) = 1_B$, since $u$ is unitary.
\end{proof}

\begin{lemma} \label{oacas}  Let $T :  A \rightarrow B$
be a linear map between unital $\oas$.   The following are
equivalent:
\begin{itemize}
\item [(i)]  $T$ is a Shilov isometry,
\item [(ii)]  $T = u \pi(\cdot)$ for a
completely isometric homomorphism $\pi$ on $A$ which is the
restriction of 1-1  *-homomorphism $C^*_e(A) \rightarrow
C^*_e(B)$, and a partial isometry $u \in B \subset C^*_e(B)$ with
$u^* u = \pi(1_A)$,
\item [(iii)]   $T$ is a complete isometry
with  $\{ 0 \} = J \subset C = C^*_e(B)$ in Theorem \ref{ncmat}
(v).
\end{itemize}

We may take $u = T(1_A)$. If these equivalent conditions hold,
then $T(A)$ contains a unitary (resp. isometry, coisometry) of
$C^*_e(B)$ if and only if $u = T(1_A)$ is a unitary (resp.
isometry, coisometry) of $C^*_e(B)$. Indeed if $1_B \in T(A)$ then
$\pi(A) \subset B$, and also $\pi$ is unital. If in addition
$T(1_A) = 1_B$, then $T$ is a completely isometric homomorphism,
and is the restriction to $A$ of a 1-1 *-homomorphism $C^*_e(A)
\rightarrow C^*_e(B)$.
\end{lemma}

\begin{proof}  Again we leave the equivalences (i)-(iii)
as an exercise.  If $1_B \in T(A)$ then $1_B \in \tilde{T}(A)$,
where $\tilde{T}$ is the triple morphism $C^*_e(A) \rightarrow
C^*_e(B)$ extending $T$.  Applying \ref{csca} we see that $u =
T(1)$ is a unitary of $C^*_e(B)$.  Again selecting $a_0 \in A$
such that $T(a_0) = 1_B$, it follows as in the proof of \ref{csca}
that $u^* = \pi(a_0) \in \pi(A)$.   Next notice that $$u T(a_0^2)
= u^2 \pi(a_0^2) = u T(a_0) \pi(a_0) = u \pi(a_0) = T(a_0) =
1_B.$$ Thus $u^* = T(a_0^2)$, whence $\pi(A) = u^* T(A) \in B^2
\subset B$. The other assertions are left as an exercise.
\end{proof}

If we drop the hypothesis that $T$ is a Shilov isometry, then
things become less simple.   The problem is that our abstract
result \ref{ncmat} is difficult to use if the `Shilov boundaries'
are complicated, as is the case with $H^\infty(\mathbb{D})$ say.
There is however a class of complete isometries which are not
Shilov isometries, but which are fairly tractable. Namely,
following Matheson \cite{Math}, we say that a linear complete
isometry $T : A \rightarrow B$ between unital operator algebras,
is {\em left type 1}, if $B \cap J = \{ 0 \}$ for 
the left ideal $J$ in Theorem \ref{ncmat}
(v).   Equivalently, $T$ is left type 1 if
 the restriction to $B$ of the canonical
quotient map $C^*_e(B) \rightarrow C^*_e(B)/J$, is 1-1;
 or if the map $b \mapsto b (1-p)$ on $B$ is 1-1, where $p$ is
as in the proof of \ref{ncmat}.

We remark that if one replaced $C^*_e(B)$ by any $\ca$ containing
$B$, most of our proofs below will still go through. In the
corollaries below we will however regard $B$ as a unital
subalgebra of $C^*_e(B)$, and any adjoints $b^*$ of elements in
$B$ are taken in that $\ca$.

There is an analogous definition of {\em right type 1}, namely
that $B \cap (NB)  = \{ 0 \}$, in the language of the proof of
\ref{ncmat}.   One can show that $T$ is  right type 1 if and only
if the map $T^* : A^* \rightarrow B^*$ given by $T^*(a^*) =
T(a)^*$, is left type 1.

We say that $T$ is {\em type 1} if it is both left and right type
1.

Note that if the function algebras in \ref{mat} are $C(K)$ spaces,
then `most' isometries are type 2 (non-type 1) isometries (since
the
type 1 isometries between $C(K)$ spaces are exactly the maps characterized
in \ref{csca}).  Quite obviously, any  Shilov isometry between
function or operator algebras is type 1, since
in this case $J = \{ 0 \}$.   However it is easy to find type 1
isometries that are not Shilov isometries.  For example, Alec
Matheson showed us an interesting class of examples constructed
along the following lines. Let $A = A(\mathbb{D})$ be the disk
algebra, and let $\psi$ be the Riemann map from the closed disk
onto the upper half of the closed disk.  Then let $\varphi(z) =
\psi(z)^2$, and define $Tf(z) = f(\varphi(z))$ for $f \in A$ and
$z$ in the closed unit disk.   It is easy to see that $T : A
\rightarrow A$ is a unital type 1 isometry (it is type 1 because
of the well known fact that a nontrivial function in
$A(\mathbb{D})$ cannot vanish on a nontrivial arc of the circle),
which is not a Shilov isometry.

For us, the main interest in type 1 isometries is the following point: 
Ideally one would like to classify isometries in terms of
composition operators. But the question of what we mean by the
term `composition' needs to be clarified.    For example,
classical $H^\infty(\mathbb{D})$ lives inside $C(\partial
H^\infty) = L^\infty(\mathbb{T})$ and hence if in this context we
want to describe isometries in terms of composition operators, we
need to specify composition by what? For a result to earn the
title of a true $H^\infty$ result, such composition should be
described in terms of point transformations on $\mathbb{D}$ rather
than $\partial H^\infty$ or $\mathbb{T}$. The papers
\cite{Math,EW} suggest that for common function algebras (on
domains in $\mathbb{C}$ for example), it is the type 1 isometries
between these algebras which seem to have 
some hope of being classifiable as composition
operators on their domains (as opposed to on their Shilov
boundaries).   What makes this work in \cite{Math} is the
commutative  analogue of our next result (see also the remarks at the 
end of this  
section for a demonstration of this technique of Matheson in a
concrete case).   Thus our next result should be useful in the
study of particular noncommutative operator algebras.

\begin{corollary}  \label{ncmatco2}  Let $T :  A \rightarrow B$
be a (not necessarily surjective) left type 1 linear complete
isometry between unital operator algebras such that $T(1)$
commutes with $T(a)$ for all $a \in A$. Then $$T(1) T(ab) = T(a)
T(b)$$ for all $a, b \in A$. Therefore if, further, $T(1)$ is
invertible, then we may write $$T = T(1) \theta(\cdot)$$ for a
unital 1-1 `completely bicontinuous' homomorphism $\theta$ defined
on $A$, namely $\theta = T(1)^{-1} T(\cdot)$.
\end{corollary}

\begin{proof}    It suffices to show that
$T(1) T(ab) (1-p) = T(a) T(b) (1-p)$. We obtain using (*) (see the
comment following \ref{ncmat}) that $$T(1) T(ab) (1-p) = T(1) T(a)
T(1)^* T(b) (1-p).  $$ Again using (*), and the fact that $T(1)$
commutes with $T(a)$, we get $$T(1) T(ab) (1-p) = T(a) T(1) T(1)^*
T(b) (1-p) = T(a) T(b) (1-p).$$

To see the last assertion note that in the noncommutative case,
$$T(1)^{-1} T(a) T(1)^{-1} T(b) = T(1)^{-2} T(a) T(1) T(1)^{-1}
T(b) = T(1)^{-1} T(ab)$$ since $T(ab) = T(1)^{-1} T(a) T(b)$.
\end{proof}

\begin{corollary}  \label{ncmatco3}  Let $T :  A \rightarrow B$
be a left type 1 unital linear complete isometry between unital
operator algebras.   Then $T$ is a homomorphism.
\end{corollary}

\begin{corollary}  \label{ncmatg}  Let $T :  A \rightarrow B$
be a left type 1 linear complete isometry between unital $\oas$.
Suppose that $T(1_A) \in \Delta(B)$.  Then $u = T(1_A)$ is a
partial isometry in $C^*_e(B)$, and there exists a completely
isometric homomorphism $\theta : A \rightarrow B$, such that $$T =
u \theta(\cdot) \quad \mbox{and} \quad \theta = u^* T(\cdot).$$
\end{corollary}

\begin{proof}    If we define $\theta = u^* T(\cdot)$, then $\theta :
A \rightarrow B$. Also, $$u \theta(a)(1-p) =  u u^* T(a) (1-p) =
T(a) (1-p)$$ by (*).  Since we are in the left type 1 situation
the last equation implies that $T = u \theta(\cdot)$. Similar
considerations show that $\theta$ is a homomorphism. Similarly we
may also conclude from (*) and the left type 1 hypothesis that $u
u^* u = u$, so that $u$ is a partial isometry in $C^*_e(B)$. Then
$\theta$ is a complete isometry, since for example $$\Vert a \Vert
= \Vert T(a) \Vert \geq \Vert u^*T(a) \Vert = \Vert \theta(a)
\Vert \geq \Vert u \theta(a) \Vert = \Vert T(a) \Vert = \Vert a
\Vert.$$
\end{proof}

\begin{corollary}  \label{ncmato}  Let $T :  A \rightarrow B$
be a left type 1 linear complete isometry between unital $\oas$.
Suppose that $1_B \in $ Ran $T$.  Then $u = T(1_A)$ is a
coisometry with $u, u^* \in B$ (that is, $u \in \Delta(B)$), and
there exists a completely isometric homomorphism $\theta : A
\rightarrow B$, possessing all the properties of the previous
corollary. Also $u^* \in \theta(A)$. If, further, $T(1_A)$
commutes with $T(A)$ (or even simply with $T(a_0^2)$, where
$T(a_0) = 1_B$), or if $T$ is also right type 1,  then $u = T(1)$
is a unitary of $\Delta(B)$, and $\theta(1_A) = 1_B$.
\end{corollary}

\begin{proof}    If $T(a_0) = 1_B$
then we have $\Vert a_0 \Vert = 1$.  As in the first displayed
equation in the proof of \ref{ncmatco2} (with $a = b = a_0$), we
have from (*) that $$u T(a_0^2) (1-p) = u  T(a_0) u^* T(a_0) (1-p)
= u u^* T(a_0) (1-p) = T(a_0) (1-p).$$ Thus $T(1) T(a_0^2) =
T(a_0) = 1_B$.  Since $T(1)$ and $T(a_0^2)$ are contractions on
some Hilbert space, this forces $u = T(1)$ to be a coisometry, and
$T(a_0^2) = T(1)^*$.   Thus $T(1)^* \in B$, and we may apply
\ref{ncmatg}. We have $\theta(a_0) = u^* T(a_0) = u^*$.

If $u = T(1)$ commutes with $T(a_0^2) = u^*$, then it is clear
that $u$ is a unitary. Similarly, if $T$ is both `left-' and
`right type 1', then a symmetric argument to the above shows that
$u$ is an isometry: $$(1-q) T(a^2_0) T(1) = (1-q) T(a_0) u^*
T(a_0) u = (1-q) T(a_0) u^*  u = (1-q) T(a_0) = (1-q).$$ Thus
$T(a^2_0) T(1) = 1_B$, so that $u$ is an isometry,  and $u$ is
unitary.
\end{proof}

Returning to Corollary \ref{ncmatco2}, we note that if in addition to the
 hypotheses there,
$B$ is a commutative algebra and
if $T(1)$ is invertible in $B$, then  the homomorphism $\theta$ in
that corollary maps into $B$, so that $A$ is commutative too.   If
in addition, $B$ is a function algebra, then so is $A$, and from
Gelfand theory it then follows that the homomorphism $\theta$ in
that corollary is a contraction.  Consequently it is an
isometry, since $$\Vert \theta(a) \Vert \geq \Vert T(1) \theta(a)
\Vert = \Vert T(a) \Vert = \Vert a \Vert.$$

\begin{corollary}  \label{ncmatco4}  Let $T :
A \rightarrow B$  be a left type 1 linear complete isometry
between unital operator algebras, and suppose that $B$ is
commutative. If $T(A)$ does not vanish identically at any point of
the maximal ideal space of $B$, then $A$ is commutative, $T(1)$ is
invertible in $B$ and $$T \; = \; T(1) \; \theta(\cdot)$$ for a
unital completely bicontinuous
 homomorphism $\theta : A \rightarrow B$.
\end{corollary}

\begin{proof}    If $\chi$ is a character of $B$ with
$\chi(T(1)) = 0$, then by \ref{ncmatco2} it follows that
$\chi(T(a)^2) = \chi(T(1)) \chi(T(a^2)) = 0$ for all $a \in A$.
Thus $\chi(T(A)) = 0$, which is a contradiction. By Gelfand
theory,  $u$ is invertible in $B$.   The rest follows from the
remark before the corollary, and \ref{ncmatco2}.
\end{proof}

\begin{corollary}  \label{ncmatco5}  Let $T :
A \rightarrow B$  be a type 1 linear isometry between
unital function algebras, and suppose that $B$ is a closed unital
subalgebra of $C(K)$ for a compact Hausdorff space $K$.  If $T(A)$
does not vanish identically at any point of $K$,  then $T(1)$ is
nonzero at every point in $K$ and $$T \; = \; T(1) \;
\theta(\cdot)$$ for a unital isometric homomorphism $\theta : A
\rightarrow C(K)$.
\end{corollary}

\begin{proof}    The argument in the previous corollary
shows that $u = T(1)$ is nonzero on $K$, and hence it is
invertible in $C(K)$.   The rest follows from the remark
 preceding \ref{ncmatco4}.
\end{proof}

There are surely noncommutative versions of the previous two
corollaries. One such result goes as follows:

\begin{corollary}  \label{ncmatco6}  Let $T :  A \rightarrow B$
be a (not necessarily surjective) left type 1 linear complete
isometry between unital operator algebras such that $T(1)$
commutes with $T(a)$ for all $a \in A$. Suppose that $B$ is a
unital subalgebra of $B(H)$, suppose that $T(A) H$ is dense in
$H$, and suppose that there is a constant $K > 0$ such that for
all $\zeta \in H$ $$\Vert \zeta \Vert \leq K \sup \{ \Vert T(a)
\zeta \Vert : a \in Ball(A) \}.$$ Then $T(1)$ is an invertible
operator on $H$, and we may write $$T = T(1) \; \theta(\cdot)$$
for a unital `completely bicontinuous' homomorphism $\theta : A
\rightarrow B(H)$. If indeed $K = 1$, then $T(1)$ is unitary and
$\theta$ is completely isometric.
\end{corollary}

\begin{proof}  By \ref{ncmatco2}, we need only show that $T(1)$ is
invertible.   We have $$\Vert \zeta \Vert \; \leq \; K^2 \sup_{b
\in Ball(A)} \sup_{a \in Ball(A)} \{ \Vert T(a) T(b)  \zeta \Vert
\} \leq K^2 \Vert T(1) \zeta \Vert$$ using the fact that $T(a)
T(b) = T(ab) T(1)$ from \ref{ncmatco2}.  Thus $T(1)$ is
bicontinuous. and has closed range.     However
 $$T(A) H = T(A) \overline{T(A) H} \subset
\overline{T(A) T(A) H} \subset \overline{T(1) H}$$ (by continuity
and the fact that $T(a) T(b) = T(1) T(ab)$). So $T(1)$ maps onto
$H$. Thus $T(1)$ is invertible.
\end{proof}


We end this section and Part I, by describing very briefly some results 
found together 
with Alec Matheson.  These 
are applications of the ideas  
above, in an attempt to gain insight into  the difficult open problems
concerning isometries $T : A \rightarrow A$, where $A =
H^\infty(\mathbb{D})$, and $\mathbb{D}$ is the open unit disk.
We first remark that it is easy to characterize 
contractive homomorphisms of $H^\infty(\mathbb{D})$, by the methods 
of p. 144 in \cite{Ho}.  These are exactly the 
composition operators $T(f) = f \circ \tau$ on $\mathbb{D}$,
for a $\tau \in Ball(H^\infty(\mathbb{D}))$.  If $T$ is 
also an isometry then one can say quite a bit more about $\tau$.
Then using the ideas above, one can show that a general type 1 isometry
$T : H^\infty(\mathbb{D})  \rightarrow H^\infty(\mathbb{D})$ 
is of the form $T(f)  = T(1) \; f \circ \tau$ on $\mathbb{D}$, 
for all $f \in H^\infty(\mathbb{D})$, and where 
$\tau$ is as above (i.e. the map $f \mapsto f \circ \tau$ is an isometric
homomorphism on $H^\infty(\mathbb{D})$).

Leaving the type 1 case, and considering general unital isometries 
$T : H^\infty(\mathbb{D})  \rightarrow H^\infty(\mathbb{D})$, 
we note firstly that because the injective envelope and 
triple envelope of $H^\infty(\mathbb{D})$ coincide with 
$L^\infty(\mathbb{T})$ (see the Remark before \ref{wsd} below), 
any such $T$ extends to a unital isometry
$\tilde{T} : L^\infty(\mathbb{T}) \rightarrow L^\infty(\mathbb{T})$.
We obtain from \ref{mat} a closed ideal $I$ in $L^\infty(\mathbb{T})$,
as in that theorem.  We let $J = H^\infty(\mathbb{D}) \cap I$.  
Denoting by $q_I$ and $q_J$ respectively the canonical quotient
maps $L^\infty(\mathbb{T}) \rightarrow L^\infty(\mathbb{T})/I$ 
and $H^\infty(\mathbb{D}) \rightarrow H^\infty(\mathbb{D})/J$,
we have that $q_J \circ T$ is an  isometric homomorphism, and 
$q_I \circ \tilde{T}$ is a 1-1 *-homomorphism (all this in fact 
is true in much greater generality).  
Let $\imath$ and $\tau$ respectively denote
the functions $\imath: z \rightarrow z$ and $T(\imath)$ on
$\mathbb{D}$.   Then $\tau = T(\imath)$ is a self-map on $\mathbb{D}$ which
induces a composition operator $C_{\tau}: H^\infty \rightarrow
H^\infty : f \mapsto f \circ \tau$, which approximates the action
of $T$ modulo the ideal ${J}$ in many senses.    For example,
$q_J \circ C_\tau$ agrees with $q_J \circ T$ on the disk algebra.
Then if $\{\delta_z: z \in \mathbb{D}, \delta_z \in
{J}^\perp\}$ is weak* dense in $\frak{M}(H^\infty) \cap {J}^\perp$
(where $\frak{M}(H^\infty)$ denotes the maximal ideal space of
$H^\infty$), then the operators $q_J \circ C_\tau$ and $q_J \circ T$
agree on all of $H^\infty$.  In this connection we remark that 
the famous Carleson Corona Theorem \cite{Car} guarantees the weak* 
density of $\{\delta_z: z\in \mathbb{D}\}$ in the maximal ideal space 
$\frak{M}(H^\infty)$ of $H^\infty$, so that our condition seems 
interesting.   Indeed Carleson's theorem in this form gives another 
route to the characterization of unital type 1 isometries mentioned in the 
previous paragraph.     This is work in 
progress, and details will be forthcoming.   

\vspace{4 mm}
\begin{center}
{\Large Part II.   Logmodularity and the $C^*$-envelope.}
\end{center}
\section{Logmodularity and representing measures.}
A Dirichlet algebra is a unital subalgebra of $C(K)$ such that $A
+ \bar{A}$ is norm dense in $C(K)$ where here $\bar{A}$ denotes
the set of adjoints of $A$. It follows that $A$ separates points
of $K$, and that $K = \partial A$, the Shilov boundary of $A$. We
therefore  define a (noncommutative) {\em  Dirichlet algebra} to
be a unital subalgebra $A$ of a unital $C^*$-algebra $B$, such
that $Re(A)$ is norm dense in $B_{sa}$.  Notice that this is
easily seen to be equivalent to saying that $A + A^*$ is norm
dense in $B$. 

In order to define a noncommutative version of logmodularity, we
will first recall a few basic facts about $C^*$-algebras. First,
we recall that in a unital $C^*$-algebra, an element $a$ is
strictly positive if and only if it is positive and invertible,
and if and only if $b \geq \epsilon I$ for some real $\epsilon >
0$.  It will be helpful also to recall the following fact
(which may be proved  for example with the assistance
of \cite{Ped} 1.3.8):
If $b$ is a strictly positive element in a  unital $C^*$-algebra,
and if $a_n$ is a sequence  of
positive elements in the algebra,
then $a_n  \rightarrow b$ uniformly
(i.e. in the norm topology) if
and only if $\sqrt{a_n} \rightarrow \sqrt{b}$.

We now proceed to define several related concepts, valid for
a unital subalgebra $A$ of a unital $C^*$-algebra $B$.
Firstly, we say that $A$ has {\em factorization}, if each strictly
positive $b \in B$ may be written as $a^* a$ for some $a \in
A^{-1}$.  This notion has been studied by many authors (e.g.
\cite{Pitts,Po}). Next, we say that $A$ is {\em left approximating
in modulus} (resp. {\em left convexly approximating in modulus})
if every positive $b \in B$ is a uniform limit of terms of the
form $a^* a$ for an $a \in A$ (resp. $\sum_{k=1}^n a_k^* a_k$ for
$a_k \in A$, $n$ varying too). Thus if $P = \{ a^* a : a \in A \}
\subset B_+$, then $A$ is left approximating in modulus (resp.
left convexly approximating in modulus) if and only if $\bar{P} =
B_+$ (resp. $\overline{\text{co}\{P\}} = B_+$).   The word `left'
here refers to the preference of products $a^* a$ as opposed to $a
a^*$; thus right approximating in modulus means that each positive
$b \in B$ is a uniform limit of terms of the form $a a^*$ for $a
\in A$.  If the word left or right is omitted, then we mean that
both left and right hold.

We define a (noncommutative) {\em logmodular} algebra to be a
unital subalgebra $A$ of a unital $C^*$-algebra $B$, such that
every strictly positive element $b \in B$  is a uniform limit of
terms of the form $a^* a$ where $a \in A^{-1}$.

Some well known results for algebras with `factorization' carry over to
the `logmodular case':

\begin{proposition} \label{acfl}  Let $A$ be a unital
subalgebra of a  unital
$C^*$-algebra $B$.  The following are equivalent:
\begin{itemize}
\item [(i)]  $A$ has factorization, i.e.
every strictly positive element
$b \in B$ may be factored $b = a^* a$ for some $a \in A^{-1}$,
\item [(ii)]  every strictly positive element
$b \in B$ equals $|a|$  for some $a \in A^{-1}$,
\item [(iii)]   every $b \in B^{-1}$ equals $u |a|$ for an $a \in A^{-1}$
and a unitary $u \in B$  (in fact $u$ may be taken to be the
unitary in the polar decomposition of $b$),
\item [(iv)]  every $b \in B^{-1}$ equals $u a$ for an $a \in A^{-1}$
and a unitary $u \in B$.
\end{itemize}
Consequently if $A$ has factorization, and if
$b \in B^{-1}$, then
$b A b^{-1} = u A u^*$ for a unitary $u \in B$.

Also the following are equivalent:
\begin{itemize}  \item [(i)$'$]  $A$ is logmodular,
\item [(ii)$'$]  every strictly positive element
$b \in B$ is a uniform limit of terms of the form
$|a|$ for $a \in A^{-1}$,
\item [(iii)$'$]   every $b \in B^{-1}$ is a uniform limit
of terms of the form $u |a|$, where $a \in A^{-1}$ and $u$ is a
unitary in $B$ (in fact $u$ may be taken to be the
unitary in the polar decomposition of $b$),
\item [(iv)$'$]   every $b \in B^{-1}$ is a uniform limit
of terms of the form $u a$, where $a \in A^{-1}$ and $u$ is a
 unitary in $B$.
\end{itemize}
 \end{proposition}

\begin{proof}   That (i) is equivalent to (ii) is
obvious, whereas the fact that (i)$'$ is equivalent to (ii)$'$
follows from the remarks made in the 
second paragraph of
this section.

 If $b \in B^{-1}$, then $b^* b \in B^{-1}$.  Thus
$b^* b$, and consequently $|b|$, is strictly positive. We may
polar decompose $b = u |b|$, with $u$ a unitary in $B$. Supposing
(i) to be true, we may write $b^* b = a^* a$ for $a \in A^{-1}$,
giving $b = u |a|$, and thereby establishing (iii).  Given (iii),
if $a = w |a|$ is the polar decomposition of $a$, then $b = u w^*
|a|$, giving (iv). Given (iv) and a  strictly positive $b$, then
$\sqrt{b} = u a$ for an $a \in A^{-1}$, so that $b = a^* a$ as in
(i).

The other assertions are similarly proved.
\end{proof}

Parts (i)-(iv) are essentially in \cite{Arvan,Pitts} (see
section
4.2 in \cite{Arvan} for example, where it is explained that  the
$u$ in (iv) is unique in a certain sense).

There seems to be a concept situated somewhere between
logmodularity and convexly approximating in modulus (see
\cite{BJ}): we say that a uniformly closed unital subalgebra $A$
of a unital $\ca$ $B$ is {\em logrigged} if every strictly
positive element $b \in B$  is a uniform limit of terms of the
form $\sum_{k=1}^n a_k^* a_k$, with $a_k \in A$, and $n$ varying
too, where there exists $b_k \in A$ with $\sum_{k=1}^n b_k a_k =
1$, and with the expressions $\sum_{k=1}^n b_k b_k^*$ converging
uniformly to $b^{-1}$.  By the remark in the 
second paragraph
at the start of this
section, this is equivalent to saying that every strictly positive
$b \in B$ is a uniform limit of terms of the form $(\sum_{k=1}^n
a_k^* a_k)^{\frac{1}{2}}$, and $b^{-1}$ is a uniform limit of
$(\sum_{k=1}^n b_k b_k^*)^{\frac{1}{2}}$, where $a_k, b_k \in A$
with $\sum_{k=1}^n b_k a_k = 1$ ($n$ varying too).

As yet we have not been able to construct an example of an algebra
which is logrigged but not logmodular.   Nonetheless since several of
our proofs below are not more difficult for logrigged algebras than
for logmodular algebras, we will state these results for the former
class.

There are various relationships between these notions, some of which
are trivial (e.g. `approximating in modulus'
 $\Rightarrow$ `convexly approximating in modulus').
The main implications to bear in mind are:

\begin{proposition} \label{allimp}  If $A$ is a unital
subalgebra of a unital $C^*$-algebra, then we have the
following implications:

`Factorization'  $\Rightarrow$
`logmodular' $\Rightarrow$ `logrigged' $\Rightarrow$
`convexly approximating in modulus'.
 \end{proposition}

\begin{proof}
The proof of the first three implications are trivial.
To see the last implication, one needs to
note also that if $b \in B$, $b \geq 0$, then
the terms
 $b + \frac{1}{n} 1$ are strictly positive and converge
uniformly to $b$.
\end{proof}

{\bf Remarks and examples.}   1)
For a function algebra $A$, the definitions above
coincide with the classical ones.   `Factorization' is
sometimes called `strongly logmodular'.

2)  The condition in the definition of a (noncommutative)
logmodular algebra that every strictly positive $b$ is a uniform
limit of terms of the form $a^* a$, implies that every strictly
positive $b$ is a uniform limit of terms $a a^*$.    This follows
by replacing $b$ with $b^{-1}$, and noting that the inverse of $a
a^*$, if $a$ is invertible in $A$, is $ (a^{-1})^* a^{-1}$.  Thus
the `logmodular' condition is more symmetric than it appears at
first sight, and there is no need to consider `left' or `right'
logmodular.

A similar remark holds for the `factorization', or
`logrigged', definition.

3)
The algebra of $n \times n$ upper triangular
matrices is Dirichlet, and is known to have factorization (this is
the Choleski factorization). Thus it is logmodular.  This can be
generalized to certain nest algebras (see e.g. \cite{Pitts,Po}).
The $H^\infty$ algebras of \cite{SW} and their
noncommutative generalization in \cite{Arvan} have factorization.

4)  
All logmodular (and most logrigged) algebras will furnish examples
of the strong Morita equivalence of the first author
with Muhly and Paulsen,  as may be seen
by the ideas in \cite{BJ} (see the end of Section 6 there).

\begin{proposition}  \label{lmis}  Suppose that $A$ is
a unital subalgebra of a unital $C^*$-algebra $B$, which is either
Dirichlet, or is left or right convexly approximating in modulus.
Then $B = C^*_e(A)$.
\end{proposition}

\begin{proof}    The `Dirichlet' assertion
is in \cite{Arv1}, but in any case is quick to verify: namely the
canonical *-epimorphism $B \rightarrow C^*_e(A)$ is, by \ref{Arv}, an
isometry and therefore is  1-1.

Suppose $A$ is left convexly approximating in modulus (the `right'
case will be similar).   By the Arveson-Hamana theorem
mentioned in the introduction there is a *-homomorphism $\pi$ from $B$ onto
$C^\ast_e(A)$.    Let $I$ be the kernel of $\pi$, we will show that
$I = \{ 0 \}$.  Let
$q$ be the canonical map $q : A \rightarrow B/I$, factoring through the
canonical maps $A \rightarrow B \overset{q_I}{\rightarrow} B/I$.
Since the complete isometry $\pi_{|_A} = j$ factors through
$q$, it follows that $q$ is a complete isometry.

Suppose $b \in I$ with $b \geq 0$. Then
$b$ is a limit of terms of the form $\sum_{k=1}^n a_k^* a_k$, and
hence $Q(b)$ is the limit of terms $\sum_{k=1}^n Q(a_k)^* Q(a_k) =
\sum_{k=1}^n q(a_k)^* q(a_k)$, since $Q$ is a *-homomorphism.
Also, $\Vert b \Vert$ is a limit of terms $ \Vert \sum_{k=1}^n
a_k^* a_k \Vert$.  However the last quantity is the square of the
norm of the column in the column-space $C_n(A)$ (that is,
the operator space given by the first column of $M_n(A)$)
whose $k$th entry
is $a_k$. Since $q$ is a complete isometry, this norm coincides
with the square of the norm of the column in $C_n(B/I)$ whose $k$th
entry is $q(a_k)$.  Thus $$\Vert b \Vert \; = \; \lim \Vert
\sum_{k=1}^n q(a_k)^* q(a_k) \Vert = \Vert Q(b) \Vert = 0.$$
Clearly $b = 0$, which implies that $I = \{ 0 \}$.
\end{proof}

Below we use the fact alluded to in the introduction that a unital
linear map $S$ between operator systems (or unital $C^*$-algebras)
is completely contractive if and only if it is completely
positive, in which case it is *-linear: i.e.  $S(v^*) = S(v)^*$.

\begin{theorem} \label{lmis2}  Suppose that $A$ is
either a Dirichlet or logrigged subalgebra of a unital
$C^*$-algebra $B$. Then any unital completely contractive
(resp. completely isometric)
homomorphism $\varphi : A \rightarrow B(H)$ admits of a unique
extension to a completely positive and completely
contractive (resp.
and completely isometric) map $B \rightarrow B(H)$.
More generally if $C$ is a
unital operator algebra, and if $\varphi : A \rightarrow C$ is a
unital contractive homomorphism which possesses an extension to a
completely positive map $B \rightarrow C$, then that extension is unique.
\end{theorem}

\begin{proof}  In the Dirichlet algebra case the claim regarding
the existence of an extension follows immediately from \ref{Arv}.
So assume that $A$ is logrigged. By the injectivity
of $B(H)$, there does exist a
completely positive extension to $B$.

To prove the claim regarding uniqueness suppose that $\Phi$ and
$\Psi$ are two completely positive extensions of $\varphi: A \rightarrow C$
to all of $B$. Since $C$ is a
unital subalgebra of some $B(H)$, it suffices to consider only the
case $C = B(H)$. Let $\xi$ be a unit vector in $H$, and suppose
that $a_k , b_k \in A$ with $\sum_{k=1}^n b_k a_k = 1$. Then $$1 =
\langle \xi , \xi \rangle = \langle \varphi(\sum_{k=1}^n b_k a_k)
 \xi , \xi \rangle = \sum_{k=1}^n
\langle \Phi(b_k) \Psi(a_k) \xi , \xi \rangle.$$
 Using the Cauchy-Schwarz inequality in a standard way yields
$$ 1 \leq \langle \sum_{k=1}^n \Psi(a_k)^* \Psi(a_k) \xi , \xi
\rangle \;  \langle \sum_{k=1}^n \Phi(b_k) \Phi(b_k)^*  \xi , \xi
\rangle,$$ which by the Kadison-Schwarz inequality (see
\cite{ERbook} 5.2.2) gives $$1 \; \leq \;  \langle \Psi \left(
\sum_{k=1}^n a_k^* a_k \right) \xi , \xi \rangle \; \langle \Phi
\left( \sum_{k=1}^n b_k b_k^* \right) \xi , \xi \rangle.$$ Since
$A$ is logrigged, this yields $$1 \; \leq \;  \langle \Psi(b) \xi
, \xi \rangle \; \langle \Phi(b^{-1})
 \xi , \xi \rangle$$
for all strictly positive $b \in B$.  Writing $b = e^u$ for $u \in
B_{sa}$, we may then replace $u$ with $tu$ for real $t$, to obtain
$$1 \; \leq \;  \langle \Psi(e^{tu}) \xi , \xi \rangle \;
\langle \Phi(e^{-tu}) \xi , \xi \rangle .$$
Let $f(t) = \langle \Psi(e^{tu}) \xi , \xi \rangle \;
\langle \Phi(e^{-tu}) \xi , \xi \rangle$, and
differentiate as in the classical proof (see e.g. 17.1 in
\cite{sto}).   We get $f'(0) = 0$ which gives
$\langle \Psi(u) \xi , \xi \rangle =  \langle
\Phi(u) \xi , \xi \rangle$.
 Hence $\Psi(u) = \Phi(u)$ for all self-adjoint
$u \in B$. From this the claim regarding uniqueness of extensions
is clear.

If $\varphi$ is completely isometric then so is
any extension of it to $B$ by the `essential' property of
$C^*_e(A)$ (see \cite{Ham}).
\end{proof}

{\bf Remark:} There is a version of the above result for general
unital contractions. In the Dirichlet case essentially the same
argument shows that such a map admits of a positive extension to
all of $B$. In proving the uniqueness of such extensions in the
logrigged case, note that we only needed the fact that $\Phi$ and $\Psi$
are continuous and
satisfy the Kadison-Schwarz inequality $\Phi(x)^* \Phi(x) \leq
\Phi(x^* x)$ for all $x \in A \cup A^*$. 

\begin{corollary}  Suppose that $A$ is
a Dirichlet or logrigged subalgebra of a unital $C^*$-algebra $B$.
Then (on relaxing the irreducibility requirement in \cite{Arv1})
every *-representation of $B$ is a `boundary representation' for
$A$ in the language of \cite{Arv1}.
\end{corollary}

Recall that in the commutative situation, the irreducible
boundary representations of a uniform algebra $A \subset
C(K)$, are precisely the
point evaluations $\epsilon_x$ for $x \in K$, whose restriction
to $A$ possess a
unique representing measure.  Such points $x$ comprise precisely
the Choquet boundary of $A$, and the closure of the
Choquet boundary is the Shilov boundary.
Since the `representing measures' above are precisely
completely positive maps extending the restriction of  $\epsilon_x$
to $A$, it is clear that our
last corollary may be interpreted as saying in some sense that
for a Dirichlet or logrigged algebra, the noncommutative
Choquet boundary
equals the Shilov boundary.

\vspace{3 mm}

 Suppose that $A$ is a unital subalgebra of a  unital $C^*$-algebra
$B$.  Then if $x \in B^{-1}$, the algebra $x A x^{-1}$ is
also a unital subalgebra of $B$, called a `similarity of $A$'.
Similarity obviously does not arise in the commutative case,
however it is very natural in the noncommutative case.

One question which seems
interesting, is the following:
if $B$ is the `noncommutative Shilov boundary' of $A$ (i.e.
$C^*_e(A) = B$), and if $x$ is invertible in $B$,
then what is the `noncommutative Shilov boundary' of
$x A x^{-1}$?  Note that the $C^*$-subalgebra of $B$ generated by
$x A x^{-1}$ may not be $B$.   However in a special case there is
a nice answer:

\begin{proposition} If  $A$ is a unital subalgebra of a  unital $C^*$-algebra
$B$, and if  $x \in B^{-1}$, then $x A x^{-1}$ has factorization
(resp. is logmodular) if and only if $A$ has factorization (resp.
is logmodular). If $A$  has factorization or is logmodular, then
$C^*_e(x A x^{-1}) = B$ for all $x \in B^{-1}$.
\end{proposition}

\begin{proof}  Suppose that $A$ has
factorization, and $x \in B^{-1}$.  By Lemma \ref{acfl}, $x A
x^{-1} = u A u^*$ for a unitary $u \in B$.  If $b$ is a strictly
positive element of $B$, then $u^* b u$ is strictly positive, and
so equals $a^* a$ for $a \in A^{-1}$. Thus  $b = u a^* a u^* = (u
a u^*)^* (u a u^*)$, so that $x A x^{-1} = u A u^*$  has
factorization.

Next suppose that $A$ is logmodular.  In this case we know that $x
= \lim_n u_n a_n$ where $u_n$ is unitary in $B$ and $a_n \in
A^{-1}$; and it follows that $\lim_n a_n^{-1} u_n^* = x^{-1}$, and
that $\{\Vert a_n^{-1} \Vert \}$ is bounded above. If $b$ is a
strictly positive in $B$ then for any fixed $n \in \mathbb{N}$ we
have $u_n^* b u_n = \lim_m (a^n_m)^* a^n_m$, for some $\{a^n_m\}_m
\subset A^{-1}$, so that $b = \lim_m (u_n a^n_m u_n^*)^* (u_n
a^n_m u_n^*)$. Since $\|u_n^\ast b u_n\| = \|b\|$ for every $n$,
we may select the $a_m^n$'s so that the entire collection $\{
a^n_m \}_{n,m}$ is uniformly bounded. Thus if $c^n_m = a_n^{-1}
a^n_m a_n \in A^{-1}$ then the collection $\{ c^n_m \}_{n,m}$ is
uniformly bounded, and it follows that there exist constants $K_1$
and $K_2$ such that $$\Vert u_n a_n c^n_m a_n^{-1} u_n^* -
 x c^n_m x^{-1} \Vert \; \leq \; K_1 \Vert u_n a_n - x \Vert
+ K_2 \Vert a_n^{-1} u_n^* - x^{-1} \Vert$$
and so by a triangle inequality argument, that
 $$b
= \lim_n \lim_m (u_n a_n c^n_m a_n^{-1} u_n^*)^* (u_n a_n c^n_m
a_n^{-1} u_n^*) = \lim_n \lim_m (x c^n_m x^{-1})^* (x c^n_m
x^{-1})$$ where $c^n_m \in A^{-1}$.  Thus $x A x^{-1}$ is
logmodular.

The assertion about the $C^*$-envelope follows immediately from
\ref{lmis}.
\end{proof}

\section{Conditional expectations and 
noncommutative $H^\infty$ spaces.}

We begin this section with a result on 
conditional expectations which we think is new.
It is a well known result due to Tomiyama (see \cite{Tom})
 that a unital
contractive projection from a unital $\ca$ $A$ onto a unital
subalgebra $B$ is completely contractive, and moreover is a
`conditional expectation' in the sense that $$P(b_1 a b_2) = b_1
P(a) b_2$$ for all $a \in A, b_1,  b_2 \in B$.  

\begin{proposition} A completely contractive unital
projection of a unital operator algebra $A$
 onto a unital subalgebra $B$ is a `conditional
expectation' in the sense that $$P(b_1 a b_2) = b_1 P(a) b_2$$ for
all $a \in A, b_1,  b_2 \in B$.
\end{proposition}

\begin{proof}  Let $P : A \rightarrow B$ be the projection.
Let $i : B \rightarrow A$  be the inclusion.    We will use the
basic properties of the injective envelope, as may be found in
\cite{Ham,ERbook} say. Let $(I(A),J)$ and $(I(B),j)$ be the
injective envelopes of $A$ and $B$ respectively. These are unital
$\cas$, and $J, j$ are unital completely isometric homomorphisms.
We may extend $j \circ P \circ J^{-1}$ to a completely contractive
unital map  $\tilde{P} : I(A) \rightarrow I(B)$ with $\tilde{P}
\circ J = j \circ P$. We may also extend $J  \circ i \circ j^{-1}$
to a completely contractive unital  map $\tilde{i} : I(B)
\rightarrow I(A)$, with $\tilde{i} \circ j = J \circ i$. Thus
$\tilde{P} \circ \tilde{i} = Id$ on $j(B)$, and hence by the
rigidity property of the injective envelope,  $\tilde{P} \circ
\tilde{i} = Id$ on $I(B)$.   Thus $Q = \tilde{i} \circ \tilde{P}$
is a unital completely contractive projection from $I(A)$ onto a
subspace of $I(A)$ which is completely order isomorphic to the
$\ca$ $I(B)$.  We have for $b \in B, a \in A$
that 
$$Q(J(a)) = \tilde{i}( \tilde{P}(J(a))) =
\tilde{i}(j(P(a)) = J(P(a)),$$ and thus
$$J(P(ba)) = Q(J(ba)) = Q(J(b)
J(a))  =
 Q(J(b) Q(J(a))),$$
where the last step 
uses a well known lemma of Choi-Effros (see the
proof of 6.1.2 in \cite{ERbook}).  
It follows that 
$$J(P(ba)) =
Q(J(bP(a))) = J(P(bP(a))) = J(b P(a)),$$ from which the result is
clear.
\end{proof}

{\bf Remark:}  The above proof provides an extension of $P$ to a
completely positive surjective map $\tilde{P} : I(A) \rightarrow
I(B)$.   One can say a little more.  We use the 
notation of the proof above.    First, if
$C^*(J(B))$ is the $\ca$ generated by $J(B)$ inside
$C^*_e(A) \subset I(A)$, then $\tilde{P}$ is
a *-homomorphism  from 
$C^*(J(B))$ onto $C^*_e(B)$.
To see this note that  $\tilde{i}$ is a complete order
isomorphism from $I(B)$ onto Ran $Q$.   Hence $\tilde{i}$ is a
*-isomorphism.   That is, $\tilde{i}(j(b_1) j(b_2)^* \cdots j(b_n))
= Q(J(b_1) J(b_2)^* \cdots J(b_n))$, or in other words,
$$\tilde{P}(J(b_1) J(b_2)^* \cdots J(b_n)) = j(b_1) j(b_2)^*
\cdots j(b_n)$$ for $b_1, \cdots, b_n \in B$.   Thus indeed
$\tilde{P}$ is a *-homomorphism $C^*(J(B)) \rightarrow C^*_e(B)$,
and $\tilde{P}(C^*(J(B))) = C^*_e(B)$.  Therefore also
$Q(C^*(J(B))) = \tilde{i}(C^*_e(B))$. It
follows from a well known  lemma of Choi that $\tilde{P}(a_1 x
a_2) = \tilde{P}(a_1) \tilde{P}(x) \tilde{P}(a_2)$, for $x \in
I(A)$ and $a_1, a_2 \in C^*(J(B))$.

We also remark that almost all of the above is true with the same 
proof, even
if $B = P(A)$ is not a subalgebra of $A$.  The 
only change is that we must amend the displayed equation in
the proposition to read
$$P(a_1 P(a_2)) = P(P(a_1) a_2) = P(P(a_1) P(a_2))$$
for all $a_1, a_2 \in A$.

\vspace{3 mm}

The following result was originally stated by Le Merdy \cite{LM1}
(see also 3.4 and 3.5 in \cite{LM2}) with the additional
hypothesis of `separate weak* continuity' of the product.  This
hypothesis was removed by the first author in \cite{BMDO}. However
in fact the last proposition gives a much simpler way to remove
the hypothesis:

\begin{corollary} (Le Merdy-Blecher)  If $A$ is a
unital operator algebra $A$ with an operator space predual, then
$A$ is completely isometrically isomorphic, via a homomorphism
which is also a homeomorphism for the weak* topologies, to a
$\sigma$-weakly closed unital subalgebra of $B(H)$, for some
Hilbert space $H$.
\end{corollary}

\begin{proof}
The dual of the canonical map $i : A_* \rightarrow A^{*}$
dualizes to a weak* continuous completely contractive unital
projection $P : A^{**} \rightarrow A$.   Let $J = $ Ker$ \; P$, a
weak* closed subspace of $A^{**}$.   Then $P$ is a conditional
expectation (by the proposition), so that $J$ is an
$\hat{A}-\hat{A}$-subbimodule of $A^{**}$.    Since the product on
$A^{**}$ is separately weak* continuous, $J$ is a 2-sided ideal of
$A^{**}$.     Also, if $F, G \in A^{**}$ then $F - P(F) \in J$, so
that $FG - P(F)G \in J$.   Thus $P(FG) = P(P(F)G) = P(F) P(G)$.
Thus $P$ is a homomorphism.   By elementary Banach space duality
principles we obtain a  completely isometric unital surjective
weak* continuous homomorphism from $A^{**}/J$ onto $A$.   By the
Krein-Smulian theorem, this homomorphism is a homeomorphism for
the weak* topologies (which implies that the multiplication on $A$
is separately weak* continuous). The result now follows from Le
Merdy's earlier result.
\end{proof}

We now discuss a situation related to logmodularity/factorization
in which conditional expectations naturally arise.  
Let $A$ be a unital subalgebra of a unital $C^\ast$-algebra $B$.
Using classical $H^\infty(\mathbb{D})$ as a model,
we may take $\Delta(A) = A \cap A^*$ to be a non-commutative
analogue of the complex scalar field. A comparison of this context
with the commutative setting (see for example \cite{Gam}, chapter
IV) suggests that at least one possible approach to a
non-commutative theory of Hardy spaces would be in terms of some
fixed homomorphism $\varphi:A \rightarrow \Delta(A) = A \cap A^*$
which is also a projection of $A$ onto $\Delta(A)$.  An extension
of such a homomorphism to a positive projection from all of $B$
onto $\Delta(A)$ may then be regarded as some sort of
non-commutative representing measure of $\varphi$. The questions
of existence and uniqueness of such ``representing measures'' and
the possible role of logmodularity in ensuring these
eventualities, immediately present themselves.
 
Let $M$ be a von Neumann algebra possessing a
faithful normal tracial state $\tau$ (which implies that $M$ is a
`finite von Neumann algebra').  
We say that $\tau$ {\em preserves} a 
map $\Phi : S \subset M \rightarrow M$, if $\tau \circ \Phi =
\tau$ on the domain $S$ of $\Phi$.   Suppose that  $A$ is a weak*
closed unital subalgebra $A$ of $M$.  Then $\Delta(A) = A \cap
A^*$ is a von Neumann subalgebra, and it is known (\cite{Tak};
V.2.36) that there exists a faithful normal conditional
expectation $\Phi$ of $M$ onto $\Delta(A)$ which is preserved by
$\tau$.

\begin{lemma}  \label{ne}
For a unital $*$-algebra $M$ with a faithful state $\tau$ and
a unital *-subalgebra $N$, there is at most one unital
`conditional expectation' from any
unital $N$-invariant subset $S$ of $M$ containing $N$,  onto
$N$, which is preserved by $\tau$.
\end{lemma}

 By `$N$-invariant' we mean that
$N S N \subset S$, and by `conditional expectation' here we
mean that $\Phi(a s b) = a \Phi(s) b$ for all $a,b \in N, s \in
S$. Since $N$ is unital, the fact that $\Phi$ is such a
conditional expectation can be shown to imply that $\Phi$
preserves the identity and also that $\Phi \circ \Phi = \Phi$.

\begin{proof}  Suppose that $\Psi$ was another
`conditional expectation' of $S$ onto $N$ which is preserved by $\tau$.
Then using the conditional expectation property we have
for $a \in S$ that
\begin{eqnarray*}
\tau(|\Phi(a) - \Psi(a)|^2) &=& \tau(\Phi(a)^* \Phi(a)) - \tau(\Phi(a)^*
 \Psi(a)) - \tau(\Psi(a)^* \Phi(a)) +
\tau(\Psi(a)^* \Psi(a)) \\
&=&
 \tau(\Phi(\Phi(a)^* a)) -
\tau(\Psi(\Phi(a)^* a)) - \tau(\Phi(\Psi(a)^*a)) +
\tau(\Psi(\Psi(a)^* a)) \\
&=& \tau(\Phi(a)^* a) - \tau(\Phi(a)^* a)
-  \tau(\Psi(a)^*a) + \tau(\Psi(a)^*a) \\
&=& 0.
\end{eqnarray*}
Since $\tau$ is faithful this  shows that $\Psi = \Phi$ on $S$.
\end{proof}

We turn now to Arveson's remarkable  
noncommutative generalization of the $H^\infty$ spaces. 
Let $M$ be a von Neumann algebra $M$
with a faithful normal tracial state $\tau$, and
suppose that  $A$ is
a weak* closed unital subalgebra $A$ of $M$.  Then
$\Delta(A) = A \cap A^*$ is a von Neumann subalgebra, and it is
known (\cite{Tak}; V.2.36) that there exists a faithful normal
conditional expectation $\Phi$ of $M$ onto $\Delta(A)$ which is
preserved by $\tau$.
We say that $A$ is a {\em finite maximal subdiagonal 
algebra}\footnote{The simplified form of the definition is due 
to Exel.},
or {\em noncommutative  $H^\infty$ space},
if further, $A + A^*$ is weak*-dense in $M$ and $\Phi$ is
multiplicative on $A$.    So in the philosophy of the above
discussion, $\Phi$ is in effect a representing measure of its
restriction to $A$.

Many examples of finite maximal subdiagonal algebra are given in
\cite{Arvan,MMS,MW} for example.   The weak* Dirichlet algebras
of \cite{SW} may be shown to all be
 finite maximal subdiagonal algebras (one may show that
in this case $\Delta(A) = \mathbb{C}1$).
If $A$ is a finite maximal subdiagonal algebra, then so is
$M_n(A)$.   In \cite{Arvan} it is shown that any finite maximal
subdiagonal algebra has factorization.  
It is consequently
logmodular, and we may apply the results proved in Section 4 
above.  From this and \ref{lmis} we may deduce the following
generalization of a well known fact for classical
$H^\infty(\mathbb{T})$:

\begin{corollary}  \label{cehi}  If
$A$ is a finite maximal subdiagonal algebra in $M$,
then $M = C^*_e(A)$, that is
$M$ is the `noncommutative Shilov boundary' of $A$.
\end{corollary}

{\bf Remark:}  It is well known that every commutative von Neumann algebra $M$ is an
injective Banach space, and is hence an injective operator space.
In this case we can say under the hypotheses of the previous
corollary, that $M$ is the `injective envelope' of $A$ too.  
(This follows easily from abstract principles
in \cite{Ham} or \cite{ERbook}).  
This is interesting because there are few
 
cases where the injective envelope of an operator space is
explicitly known.

\begin{definition}  \label{wsd}  Let $M$ be a von Neumann algebra with a
faithful normal tracial state $\tau$.
A {\em tracial subalgebra} of $M$ is a
weak* closed
unital subalgebra $A$ of $M$  for which
there exists a linear projection $\Psi$ from $A$ onto
$\Delta(A) = A \cap A^*$ which is
also a homomorphism on $A$, such that
$\tau(a) = \tau(\Psi(a))$ for all $a \in A$.
\end{definition}

\begin{theorem}
Let $A$ be a tracial subalgebra of a von Neumann algebra $M$. Then
the map $\Psi$ in definition \ref{wsd} is unique, completely
contractive, and weak* continuous.   Indeed $\Psi$ is the
restriction to $A$  of the canonical conditional expectation of
$M$ onto $\Delta(A)$ which is preserved by $\tau$.  Furthermore,
\begin{itemize}
\item [(1).]  if  $A
+ A^*$ is weak* dense in $M$, then $A$ is a finite maximal
subdiagonal algebra in $M$. Thus $A$ will then have factorization,
and will be logmodular.
\item [(2).]  if
$A$ is logmodular or logrigged then the canonical conditional expectation
from $M$ onto $\Delta(A)$ is the only positive extension of $\Psi$ to a
map from $M$ into $C$, for any $\ca$ $C$ containing (a copy of) 
$\Delta(A)$.
\end{itemize}
\end{theorem}

\begin{proof}
The conditional expectation  from $M$  onto $\Delta(A)$ restricts
to a `conditional expectation' from $A$ onto $\Delta(A)$. Clearly
$\Psi$ is a `conditional expectation' from $A$ onto $\Delta(A)$.
The first claim then follows by \ref{ne}.

To see the statement regarding the uniqueness of positive
extension in the logrigged case, it follows from the remark
following \ref{lmis2} that we only need to verify that each
positive extension of $\Psi$ satisfies the
Kadison-Schwarz inequality on $A \cup A^*$. To this end let 
$C$ be a unital $\ca$, 
let $j : \Delta(A) \rightarrow C$ be a unital *-homomorphism,
and let $\Phi
: M \rightarrow C$ be any positive extension of $j 
\circ \Psi$ to all of $M$.   For simplicity the reader may
want to take 
$C = M$ and $j = Id$ in the following.
By taking
adjoints we may conclude from $\Psi$'s action on $A$ that $\Phi$ is
unital, that it maps $A + A^*$ into $j(\Delta(A))$, satisfies the
equality $\Phi \circ j^{-1} \circ \Phi|_{A+A^*} = \Phi|_{A+A^*}$, and also acts
homomorphically on both $A$ and $A^*$. Therefore on expanding the
term $(x^* - j^{-1} \Phi(x^*))(x - j^{-1} \Phi(x)) \geq 0$ and applying $\Phi$,
it follows that $\Phi(x^*x) \geq \Phi(x^*)\Phi(x)$ for each $x \in
A \cup A^*$ as required.
\end{proof}

We close our paper with an open problem:
If $A$ is a logrigged or logmodular tracial subalgebra, when is $A^* +
A$ automatically weak* dense in the ambient von Neumann algebra?
We recall that in the 
setting of commutative weak* Dirichlet algebras considered
by Srinivasan and Wang, it is known that logmodularity and the
weak* density of $A^* + A$ in $L^\infty$ are equivalent.

We will give a result which
probably is in the right direction towards the solution of this problem.
In the commutative context of Srinivasan and Wang and other 
authors, one first proves $L^2$-density
of $A^* + A$, before combining Szego's theorem with this fact to
conclude that in fact weak* density of $A^* + A$ in $L^\infty$
pertains.  See also \cite{Lum} where this is linked very 
tightly to the principle of
\ref{lmis2}. Thus our question seems tied to the
open problems surrounding the apparent absence of a suitable
version of Jensen's inequality and Szego's theorem in the
noncommutative case (see \cite{Arvan,MW} and references contained
therein).

\begin{proposition}
Let $A$ be a logrigged tracial subalgebra of a von Neumann algebra
$M$, with 
$\Delta(A)$ contained in the center of $M$.  
Then $A + A^*$ is a dense subspace of $L^1(M, \tau)$.
\end{proposition}

\begin{proof}
Suppose, by way of contradiction,
that $A$ satisfies the hypotheses of the 
Proposition, but that
$A^* + A$ was not dense in $L^1(M, \tau)$.  By the
Hahn-Banach theorem and the duality of non-commutative
$L^p$-spaces (see e.g.
\cite{Terp}), we would be able to find some $x
\in M$ with $x \neq 0$ and $\tau(xa) = 0$ for all $a \in A^* + A$.
By taking
adjoints it is easy to see that $\tau(xa) = 0$ for all $a \in A^*
+ A$ if and only if $\tau((x + x^*)a) = 0$ and $\tau((i(x^* -
x))a) = 0$ for all $a \in A^* + A$. We may therefore assume $x$ to
be a self-adjoint element of $M$, and on suitably scaling $x$ even
that $\|x\|_\infty \leq 1$.   Then $1 + x$ is a positive element of $M$.  

If $d \in \Delta(A)$ and $a \in A^* + A$, then 
$$\tau(\Phi(xa)d) = \tau(\Phi(xad)) = \tau(xad) = 0.$$
Since this holds for all $d \in \Delta(A)$, we may conclude 
that $\Phi(xa) = 0$ for $a \in A^* + A$.    Consider the linear map 
$\Psi : a \mapsto \Phi((1+x)a)$ on $M$.   This coincides with 
$\Phi$ on $A$.  If $d \in \Delta(A), a \in M$ then using the 
conditional expectation property, the fact that $\tau$ 
preserves $\Phi$, and the facts that 
$d$ is in the center and $\tau$ is a trace, we have   
$$\tau(d^*\Phi((1+x)a^* a)d) = \tau(\Phi(d^*(1+x)a^* ad))
= \tau(d^* a (1+x) a^* d) \geq 0.$$
Since this holds for all $d \in \Delta(A)$,
it follows by elementary considerations that 
$\Phi((1+x)a^* a) \geq 0$.  Thus $\Psi$ is
positive, and therefore 
completely positive, seeing as its range is commutative \cite{P}.    
It follows from 
\ref{lmis2} that $\Phi((1+x)a) = \Phi(a)$ for all
$a \in M$.  Thus $\Phi(x M) = 0$, which implies that
$\Phi(x^2) = 0$.  Since $\Phi$ is faithful,
this  gives the contradiction  $x = x^2 = 0$.
\end{proof}

\vspace{2 mm}

Acknowledgments:   We thank Alec Matheson for drawing our
attention to the area discussed in Part I above, and for helpful discussions and
clarifications.    We also thank  the University of South Africa
for support for a  visit of the first author.

\end{document}